\documentclass[12pt,a4paper,reqno]{amsart}

\usepackage{anysize}
\usepackage[utf8]{inputenc}
\usepackage{amsmath}
\usepackage{amssymb}
\usepackage{array}
\usepackage[unicode]{hyperref}

\marginsize{2.5cm}{2.5cm}{3cm}{3cm}

\begin{document}

\title{Fast formulas for the Hurwitz values $\zeta(2,a)$ and $\zeta(3,a)$}

\author{Jesús Guillera}
\address{Department of Mathematics, University of Zaragoza, 50009 Zaragoza, SPAIN}
\email{}
\date{}

\begin{abstract}
We prove fast convergent formulas for the Hurwitz values $\zeta(2,a)$ and $\zeta(3,a)$ respectively with the help of the WZ method. In them $(a)_n$ denotes the rising factorial or Pochhammer's symbol defined by $(a)_0=1$ and $(a)_n=a(a+1)\cdots(a+n-1)$ for positive integers $n$. The Huwitz $\zeta$ function is defined by $\zeta(s,a)=\zeta(0,s,a)=\sum_{k=0}^{\infty} (k+a)^{-s}$. In addition, we can use these fast evaluations to compute also in a rapid way Dirichlet values of the kinds $L_{\chi}(2)$ and $L_{\chi}(3)$.
\end{abstract}

\maketitle

\section{Wilf-Zeilberger (WZ) pairs and continued fractions}
Herbert Wilf and Doron Zeilberger invented the concept of WZ pair: Two hypergeometric (in $n$ and $k$) terms $F(n,k)$ and $G(n,k)$ form a WZ pair if the identity
\[
F(n+1,k)-F(n,k)=G(n,k+1)-G(n,k)
\]
holds.
A Maple code written by Zeilberger, available in Maple \cite{flawless-WZ} 
finds the mate of a term (that forms a WZ pair with it), whenever such exists, by means of a rational certificate $C(n,k)$ so that $G(n,k)=C(n,k) F(n,k)$ \cite{certificate}. In \cite{flawless-WZ} we discovered WZ pairs satisfying $F(0,k)=F(+\infty,k)=0 \quad  \forall k \in \mathbb{C}$ (that we called flawless WZ pairs) and have really interesting properties. Here, we will show another kind of very interesting WZ pairs, those satisfying
\[
\lim_{k \to \infty} \sum_{n=0}^k G(n,k)=0, \quad
\lim_{n \to \infty} \sum_{k=0}^n F(n,k)=0,
\]
Then, we have
\[
\sum_{k=0}^{\infty} F(0,k+x)=\sum_{n=0}^{\infty} G(n,x),
\]
and the summation on the left is simpler (and very inteesting) if $F(0,k)$ is a rational function. The following algorithm always holds:
\[
\sum_{n=0}^{\infty} G(n,x)=G(0,x)+G(0,x)T(0,x)+G(0,x)T(0,x)T(1,x)+\cdots,
\]
where
\[
T(n,x)=\frac{G(n+1,x)}{G(n,x)}.
\]
With the help of  Euler's formula
\[
a_0 + a_0 a_1 + a_0 a_1a_2+ a_0 a_1 a_2 a_3 \cdots
=0+\frac{a_0}{1+{\displaystyle \frac{-a_1}{1+a_1+{\displaystyle
\frac{-a_2}{1+a_2+{\displaystyle \frac{-a_3}{1+a_3+\ddots}}}}}}},
\] 
we denote with $P_n(x)$ and $Q_n(x)$ the numerator and denominator of the function $T(n,x)$ after we have simplified it for a particular value of $x$. Hence
\[
T(n, x)=\frac{P_n(x)}{Q_n(x)} \quad \text{if \, \, $n \geq 0$},
\]
and
\[
A_n(x)=P_n(x)+Q_n(x), \quad B_n(x)=Q_n(x)P_{n+1}(x).
\]
Then, we see that we can write our formulas as 
\[
\sum_{n=0}^{\infty} G(n,x) = [[0,1,A_n(x)],[G(0,x),-P(0,x),-B_n(x)]] \quad \text{for \, $n \geq 0$}
\]
written in Cohen's continued fraction notation, where $n$ begins at $n=0$. A more standard notation is used in \cite{contfrac-automatic}. An excellent book about continued fractions is \cite{neverending-zudilin}. The paper \cite{Cohen-CF} includes a big database of examples. Also interesting is \cite{contfrac-automatic} where some families of continued fractions are discovered and rigorously proved in an automatic way.

\section{Simple rapid formulas for $\zeta(2, a)$ and $\zeta(3, a)$}

\subsection{A simple rapid formula for $\zeta(2, a)$}

Let
\[
U(n,k)=\frac{(1)_n^3}{\left(\frac12\right)_n (1+k)_n^2}.
\]
and
\[
S(n,k)=\frac{1}{(n+k+1)^2}, \quad R(n,k)=\frac{3n+2k+3}{2(2n+1)(n+k+1)}.
\]
Then, the pair formed by
\[
F(n,k)=U(n,k)S(n,k) \left(\frac14 \right)^n, \qquad  G(n,k)=U(n,k)R(n,k) \left( \frac14\right)^n
\]
is a WZ pair such that $F(0,k)$ is a rational function. In addition, we have
\[
\sum_{n=0}^{\infty} G(n,x)=\sum_{k=0}^{\infty} F(0,k+x)=\zeta(2, 1+x)..
\]
Explicitly, we obtain the identity
\[
\zeta(2,x)=\sum_{n=0}^{\infty} \frac{(1)_n^3}{\left( \frac12\right)_n (x)_n^2} \frac{3n+2x+1}{2(2n+1)(n+x)^2} \left( \frac14 \right)^n.
\]

\subsection{A simple rapid formula for $\zeta(3, a)$}
Let
\[
U(n,k)=\frac{(1)_n^5}{\left(\frac12\right)_n (1+k)_n^4}.
\]
and
\[
S(n,k)=\frac{n+2k+2}{2(n+k+1)^4}, \quad R(n,k)=\frac{5n^2+6nk+2k^2+10n+6k+5}{4(2n+1)(n+k+1)^4}.
\]
Then, the pair formed by
\[
F(n,k)=U(n,k)S(n,k) \left(-\frac14 \right)^n, \qquad  G(n,k)=U(n,k)R(n,k) \left(-\frac14 \right)^n
\]
is a WZ pair such that $F(0,k)$ is a rational function. In addition, we have
\[
\sum_{n=0}^{\infty} G(n,x)=\sum_{k=0}^{\infty} F(0,k+x)=\zeta(3,1+x)..
\]
Explicitly, we obtain the identity
\[
\zeta(3,x)=\sum_{n=0}^{\infty} \frac{(1)_n^5}{\left( \frac12\right)_n (x)_n^4} \frac{5n^2+6nx+2x^2+4n+2x+1}{8(2n+1)(n+x)^4} \left(-\frac14 \right)^n.
\]

\section{Accelerating our previous formulas}

\subsection{A fast convergent formula for $\zeta(2,a)$}

Applying the transformation $F(n, k+n)$ we have discovered the following WZ pair \cite{A-equal-B,certificate}:
\[
F(n,k)=U(n,k) S(n,k) \left( \frac{1}{64}\right)^n, \qquad G(n,k)=U(n,k)R(n,k) \left( \frac{1}{64}\right)^n,
\]
where 
\[
U(n,k)=\frac{(1)_n^3(1+k)_n^2}{\left( \frac12\right)_n \left(\frac12+\frac{k}{2} \right)_n^2 \left(1+\frac{k}{2}\right)_n^2}, \quad 
S(n,k)=\frac{1}{(2n+k+1)^2}.
\]
and 
\[
R(n,k)=\frac{21n^3+55n^2+47n+13+2k^3 + 13k^2n + 28kn^2 + 11k^2 + 48kn + 20k}{2(2n+k+ 1)^2(2n+1)(2n+k+2)^2}.
\]
We have 
\[
\sum_{n=0}^{\infty} G(n,x) = \sum_{k=0}^{\infty} F(0,k+x) 
= \sum_{k=0}^{\infty} \frac{1}{(k+1+x)^2}=\zeta(2,1+x),
\]
We know that $\zeta(2,1)=\zeta(2)=\pi^2/6$, and $\zeta(2,1/4)=\pi^2+8 Catalan$. Hence, we have
\[
\sum_{n=0}^{\infty} G(n,0) = \frac{\pi^2}{6}, \qquad 
\sum_{n=0}^{\infty} G(n,-3/4)=\pi^2+8 \, Catalan.
\]
Let
\[
T(n,x)=\frac{G(n+1,x)}{G(n,x)}.
\]
Simplifying $Y(n,x)$ and defining $T(n,x)$ as the explicit output of the simplification, we see that 
\begin{multline*}
T(n,x)=\frac{21n^3 + 28n^2x + 13nx^2 + 2x^3 + 118n^2 + 104nx + 24x^2 + 220n + 96x + 136}{21n^3 + 28n^2x + 13nx^2 + 2x^3 + 55n^2 + 48nx + 11x^2 + 47n + 20x + 13} \\ \times \frac{2(x + n + 1)^2(n + 1)^3}{(3 + x + 2n)^2(4 + x + 2n)^2(3 + 2n)},
\end{multline*}
and
\[
T(n,0)=\frac{(n + 1)^3(21n + 34)}{8(3 + 2n)^3(21n + 13)}.
\]
We have written the code noticing that
\[
\zeta(2, 1+x)=\sum_{n=0}^{\infty} G(n,x)= G(0,x)+G(0,x)T(0,x)+G(0,x)T(0,x)T(1,x)+\cdots,
\]
and observing the relation between a generic term and its preceding one. For $x=0$ we know that $\zeta(2,1)=\zeta(2)$.

\subsection{Maple code for computing $\zeta(2, 1+x)$ efficiently with DIG digits}

\begin{verbatim}

with(SumTools[Hypergeometric]):
Zeilberger(F(n,k),n,k,N)[1];
G:=(nn,kk)->subs({n=nn,k=kk},Zeilberger(F(n,k),n,k,N)[2]);
Y:=(n,k)->simplify(G(n+1,k)/G(n,k)):
COMPUTE:=proc(x,DIG) local t,n,H: global T,SUMA,TOTAL: t:=time(): 
# We define here T:=(n,x)-> as the explicit output of Y(n,x)#
H:=evalf(simplify(G(0,x),DIG): SUMA:=H:
Digits:=DIG: for n from 0 to floor(evalf(DIG/log(64,10))) do:
H:=H*T(n,x): SUMA:=evalf(SUMA+H,DIG): od:
print(evalf(SUMA,DIG)): print(SECONDS=time()-t): end:

\end{verbatim}

For example, for computing $\zeta(2, 1/5)$ up to 100000 digits with our formula written efficiently as the algorithm above, execute \texttt{COMPUTE(-4/5, 100000)}.

\subsection{A fast convergent formula for $\zeta(3, a)$.} Applying the transformation $F(n, k+n)$ we have discovered the following WZ pair \cite{A-equal-B,certificate}:
\[
F(n,k)=U(n,k) S(n,k) \left( -\frac{1}{1024}\right)^n, \quad G(n,k)=U(n,k)R(n,k) \left( -\frac{1}{1024}\right)^n,
\]
where 
\[
U(n,k)=\frac{(1)_n^5(1+k)_n^4}{\left( \frac12\right)_n \left(\frac12+\frac{k}{2} \right)_n^4 \left(1+\frac{k}{2}\right)_n^4}, \quad S(n,k)=\frac{3n+2k+2}{2(2n+k+1)^4},
\]
and $R(n,k)$ is obtained from the certificate. We have
\begin{align*}
\sum_{n=0}^{\infty} G(n,0) & = \frac{1}{64} \sum_{n=0}^{\infty} \frac{(1)_n^5}{\left( \frac12 \right)_n^5} \frac{205n^2+250n+77}{(2n+1)^5} \left(\frac{-1}{1024}\right)^n \\ & = \sum_{k=0}^{\infty} F(0,k) = \sum_{k=0}^{\infty} \frac{1}{(k+1)^3}=\zeta(3),
\end{align*}
due to T. Amdeberham and D. Zeilberger \cite{zeta3}. In addition, we have
\[
\sum_{n=0}^{\infty} G(n,x) = \sum_{k=0}^{\infty} F(0,k+x) 
= \sum_{k=0}^{\infty} \frac{1}{(k+1+x)^3}=\zeta(3,1+x).
\]

Acceleration with the transformation $F(n,k) \to F(n,k+n)$ leads to even a faster series.

\subsection{Maple code for computing $\zeta(3, 1+x)$ efficiently with DIG digits}

\begin{verbatim}

with(SumTools[Hypergeometric]):
Zeilberger(F(n,k),n,k,N)[1];
G:=(nn,kk)->subs({n=nn,k=kk},Zeilberger(F(n,k),n,k,N)[2]);
Y:=(n,k)->simplify(G(n+1,k)/G(n,k)):
COMPUTE:=proc(x,DIG) local t,n,H,R: global SUMA,TOTAL: t:=time(): 
# We define here T:=(n,x)-> as the explicit output of Y(n,x) #
H:=evalf(simplify(G(0,x),DIG): SUMA:=H:
Digits:=DIG: for n from 0 to floor(evalf(DIG/log(1024,10))) do:
H:=H*T(n,x): SUMA:=evalf(SUMA+H,DIG): od:
print(evalf(SUMA,DIG)): print(SECONDS=time()-t): end:

\end{verbatim}

\begin{center}
Table of times (in our home computer) for $1+x=1/5$ using our formulas \vskip 0.5cm
\begin{tabular}{|c|c| c|}
\hline 
Digits & $\zeta(2,1/5)$ & $\zeta(3,1/5)$ \\ 
\hline 
10000 & $3''$ & $2''$  \\
20000 & $42''$ & $15''$ \\
40000 & $59''$ & $33''$ \\
80000 & $292''$ & $160''$ \\
160000 & $1307''$ & $733''$ \\
\hline
\end{tabular}
\end{center}

\section{Continued fractions for $\zeta(2)$ and $\zeta(3)$}

Taking  $x=0$, we can arrive at

\[
\zeta(2)=[[0,a_{n-2}],[104,b_{n-1}]],
\]
where
\begin{align*}
a_n &=-(n + 1)^3(21n + 34) - 8(3 + 2n)^3(21n + 13), \\
b_n &=8(n + 1)^3(21n + 34)(1 + 2n)^3(21n - 8), 
\end{align*}
and
\[
\zeta(3)=[[0,a_{n-2}],[1232,b_{n-1}]],
\]
where
\begin{align*}
a_n &-(n + 1)^5(205n^2 + 660n + 532) - 32(2n + 3)^5(205n^2 + 250n + 77), \\
b_n &=32(n + 1)^5(205n^2 + 660n + 532)(1 + 2n)^5(205n^2 - 160n+32).
\end{align*}

\section{Dirichlet values $L_{\chi}(2)$ and $L_{\chi}(3)$}

The Dirichlet function $L_{\chi}(s)$ is defined as the analytic continuation of
\[
L_{\chi}(s)=\sum_{n=1}^{\infty} \left( \frac{\chi}{n} \right) \frac{1}{n^s},
\]
and is related to the Hurwitz $\zeta$ function) in the following way:
\[
L_{\chi}(s) =\frac{1}{(-\chi)^s (s-1)!} \sum_{j=1}^{|\chi|-1} \left( \frac{\chi}{j} \right) \zeta \left(s, \frac{j}{|\chi|} \right).
\]
From the last identity we observe that we can compute $L_{\chi}(2)$ and $L_{\chi}(3)$ in a fast way as well. In some cases we can even simplify the formula. For example, we have
\begin{align}
8  L_{4}(2) &= -\pi^2+ \zeta\left(2,\frac14\right), \quad
27L_{-3}(2) =-4\pi^2+6 \zeta\left(2, \frac13 \right), \\
32 L_{-8}(2) &= -4 \pi^2 + \zeta\left(2, \frac18 \right)+\zeta\left(2, \frac38 \right), 
\end{align}
and
\begin{align}
625 L_{5}(3) &= -4 \pi^3 \sqrt{25-2\sqrt{5}} + 10 \zeta\left(3, \frac15 \right)- 10\zeta \left(3, \frac25 \right), \\
512 L_{8}(3) &= -16\pi^3 +2 \zeta\left(3, \frac18 \right)-2\zeta\left(3, \frac38 \right), \\
1728 L_{12}(3) &=-32 \sqrt3 \pi^3 + 2 \zeta\left(3, \frac{1}{12} \right)-2\zeta\left(3, \frac{5}{12} \right).
\end{align}
Further useful information for the Dirichlet $L_{\chi}(s)$ is in \cite{mathworld-L}

\section{A fast formula for the Dirichlet value $L_{-8}(2)$}

There are known fast formulas suitable for computing the Dirichlet values $L_{-3}(2)$, $L_{-4}(2)$, $L_{-7}(2)$, $\zeta(3)$ and $\zeta(5)$ using $y$-cruncher. Here, we will prove another one for computing $L_{-8}(2)$. First check the following WZ pair:
\[
F(n,k)=U(n,k) S(n,k), \quad G(n,k)=U(n,k) R(n,k),
\]
where
\[
U(n,k)=\frac{(1)_n^3}{\left( \frac12+k\right)_n^3} (-1)^k (-1)^n, \qquad S(n,k)=\frac{n+2k+1}{(2n+2k+1)^3}
\]
and $R(n,k)$ is determine automatically using the package written by Zeilberger. The transformation $(n,k) \to  (n, k+n)$ leads to a WZ pair $(F, G)$, where
\[
F(n,k)=\frac{(1)_n^3 \left( \frac12+k\right)_n^3}{\left( \frac14+\frac{k}{2}\right)_n^3 \left( \frac34+\frac{k}{2}\right)_n^3} \left( \frac{1}{64}\right)^n \frac{3n+2k+1}{(4n+2k+1)^3},
\]
and $G(n,k)$ is the companion. Replacing $k$ with $k+x$ and using the properties of this WZ pair, we have
\[
\sum_{n=0}^{\infty} G(n,x) = \sum_{k=0}^{\infty} F(0, k+x)=\sum_{k=0}^{\infty} \frac{(-1)^k}{(2k+2x+1)^2}.
\]
Letting $x=1/4$, we get
\[
16 \sum_{n=0}^{\infty} G \left(n,\frac14 \right) = 16 \sum_{k=0}^{\infty} \frac{(-1)^k}{\left(2k+\frac12+1\right)^2}=\zeta \left(2, \frac38\right) - \zeta \left(2, \frac78\right).
\]
Using 
\[
32 L_{-8}(2) = -4\pi^2  + \zeta \left(2, \frac18 \right)+\zeta \left(2, \frac38 \right), \qquad 
\zeta \left(2, \frac78 \right)+\zeta \left(2, \frac18 \right) = (4+2\sqrt2) \pi^2,
\]
and making some simplifications, we arrive at
\[
L_{-8}(2)=\frac{\sqrt2}{16} \pi^2 + \sum_{n=0}^{\infty} \frac{\left(1\right)_n^3 \left(\frac34\right)_n^3}{\left(\frac38\right)_n^3 \left(\frac78\right)_n^3} \left( \frac{1}{64}\right)^n \frac{5376n^4+16768n^3+19296n^2+9660n+1761}{(8n+3)^3(8n+7)^3},
\]
which we can compute with a lot of digits. A list of fast computing series for some Dirichlet $L$-values is in \cite{zuniga-dirichlet}

\end{document}